\documentclass[a4paper, 12pt]{article}
\usepackage{amsfonts}
\usepackage{amssymb}
\usepackage{amsthm}
\usepackage[dvips]{graphics,color}

\newtheorem{theorem}{Theorem}[section]
\newtheorem{lemma}[theorem]{Lemma}
\newtheorem{defn}[theorem]{Definition}
\newtheorem{cor}[theorem]{Corollary}
\newtheorem{prop}[theorem]{Proposition}
\theoremstyle{definition}{\newtheorem{eg}[theorem]{Example}}

\newcommand{\hf}[1][1]{\frac{#1}{2}}

\newcommand{\di}{\,\mathrm{d}}
\newcommand{\rmd}{\,\mathrm{d}}
\newcommand{\Or}{\mathrm{O}}
\newcommand{\rmi}{\mathrm{i}}

\begin{document}
\setlength{\unitlength}{10mm}
\title{ A one dimensional analysis of  turbulence and its intermittence for the
$d$-dimensional stochastic Burgers equation}
\author{A D Neate $\qquad$  A Truman\\
\small{Department of Mathematics, University of Wales Swansea,}\\ [-0.8ex]
\small{Singleton Park, Swansea, SA2 8PP, Wales, UK.}}

\maketitle
\begin{abstract}

The inviscid limit of the stochastic Burgers equation is discussed in terms of the level surfaces of the
minimising Hamilton-Jacobi function, the classical mechanical
caustic and the Maxwell set and their algebraic pre-images under the classical mechanical flow map.
The problem is analysed in terms of a reduced (one dimensional)
action function. We demonstrate that the geometry of the caustic, level surfaces and Maxwell set can change infinitely rapidly causing turbulent behaviour which is stochastic in nature. The intermittence of this turbulence is demonstrated in terms of the recurrence of
two
processes.
\end{abstract}

\section{Introduction}

Burgers equation has been used in studying turbulence and in
modelling the large scale structure of the universe \cite{Arnold4,
E, Shandarin}, as well as to obtain  detailed asymptotics for
stochastic Schr\"{o}dinger and heat equations.
\cite{Elworthy,Elworthy2,Truman,Truman2,Truman3,Truman4}. It has
also played a part in Arnol'd's work on caustics and Maslov's works
in semiclassical quantum mechanics
\cite{Arnold,Arnold2,Maslov,Maslov2}.

   We consider the stochastic viscous Burgers equation for the
velocity field $v^{\mu}(x,t)\in\mathbb{R}^d$, where
$x\in\mathbb{R}^d$ and $t>0$,
\[\frac{\partial v^{\mu}}{\partial
t}+\left(v^{\mu}\cdot \nabla\right) v^{\mu} = \hf[\mu^2]\Delta
v^{\mu} - \nabla V(x) -\epsilon \nabla k_t(x)\dot{W}_t,\]
with initial condition  $v^{\mu}(x,0) = \nabla S_0 (x)+\Or(\mu^2).$
Here $\dot{W}_t$ denotes white noise
and $\mu^2$ is the coefficient of viscosity which we assume to  be small.

We are interested in the advent of discontinuities in the inviscid
limit of the Burgers fluid velocity $v^0(x,t)$  where
$v^{\mu}(x,t)\rightarrow v^0(x,t)$ as $\mu\rightarrow 0.$ Using the
Hopf-Cole transformation $v^{\mu}(x,t)=-\mu^2\nabla\ln
u^{\mu}(x,t)$, the Burgers equation becomes the Stratonovich heat
equation,
\[\frac{\partial u^{\mu}}{\partial
t}=\hf[\mu^2]\Delta u^{\mu} +\mu^{-2}V(x)u^{\mu}+\frac{\epsilon}{\mu^{2}}
k_t(x)u^{\mu}\circ\dot{W}_t ,
\]
with initial condition  $u^{\mu}(x,0) = \exp\left(-\frac{S_0(x)}{\mu^2}\right)T_0(x),$
where the convergence factor $T_0$ is related to the initial Burgers fluid
density.

Now let, \begin{equation}\label{action} A[X]: =
\hf\int_0^t \dot{X}^2(s)\di s -\int_0^t V(X(s)) \di s -\epsilon
\int_0^t k_s(X(s))\di W_s,\end{equation}
and select a path $X$ which minimises $A[X]$. This requires,
\[
\di \dot{X}(s)+\nabla V(X(s)) \di s
+\epsilon \nabla k_s(X(s))\di W_s=0.
\]
We then define the stochastic action
$A(X(0),x,t):=\inf\limits_X\left\{ A[X]:X(t)=x\right\}.$
Setting,
\[\mathcal{A}(X(0),x,t):=S_0(X(0))+A(X(0),x,t),\]
and then minimising $\mathcal{A}$  over $X(0)$, gives
$\dot{X}(0) = \nabla
S_0(X(0)).$
Moreover, it follows that, \[\mathcal{S}_t(x):=\inf\limits_{X(0)}
\left\{\mathcal{A}(X(0),x,t)\right\},\]
is  the minimal solution of the Hamilton-Jacobi
equation,
\begin{equation}\label{ie4}
\di \mathcal{S}_t +\left(\hf|\nabla \mathcal{S}_t|^2+V(x)\right)\di
t +\epsilon k_t(x) \di W_t=0, \qquad \mathcal{S}_{t=0}(x) =
S_0(x).\end{equation} Following the work of Donsker, and Freidlin
and Wentzell \cite{Freidlin}, $-\mu^2\ln u^{\mu}(x,t) \rightarrow
\mathcal{S}_t(x)$ as $\mu\rightarrow 0$. This gives the inviscid
limit of the minimal entropy solution of Burgers equation  as
$v^0(x,t)=\nabla\mathcal{S}_t(x)$ \cite{Dafermos}.

\begin{defn}
The stochastic wavefront at time $t$ is defined to be the set,
\[\mathcal{W}_t=\left\{x:\quad\mathcal{S}_t(x)=0\right\}.\]
\end{defn}
For small $\mu$ and fixed $t$,  $u^{\mu}(x,t)$ switches
continuously from being exponentially large to small as $x$ crosses the
wavefront $\mathcal{W}_{t}$. However, $u^{\mu}$ and $v^{\mu}$ can
also switch discontinuously.

Define the classical flow map
$\Phi_s:\mathbb{R}^d\rightarrow\mathbb{R}^d$ by,
\[\rmd \dot{\Phi}_s +\nabla V(\Phi_s) \rmd s+\epsilon\nabla k_s(\Phi_s)\rmd
W_s =0,\qquad\Phi_0 = \mbox{id},\qquad \dot{\Phi}_0 = \nabla S_0.\]
Since $X(t) =x$ it follows that $X(s) = \Phi_s\left(
\Phi_t^{-1}( x)\right)\!,$ where the pre-image $x_0(x,t) = \Phi_t^{-1} (x)$ is
not necessarily unique.

Given some regularity and boundedness, the global inverse function
theorem gives a caustic time $T(\omega)$ such that for
$0<t<T(\omega)$, $\Phi_t$ is a random diffeomorphism; before the
caustic time $v^0(x,t) = \dot{\Phi}_t\left(\Phi_t^{-1}(x)\right)$ is
the inviscid limit of a classical solution of the Burgers equation
with probability one.

The method of characteristics suggests that discontinuities in
$v^0(x,t)$ are associated with the non-uniqueness of the real pre-image $x_0(x,t)$.
When this occurs,  the classical flow map
$\Phi_t$ focusses
an infinitesimal volume of points $\rmd
x_0$ into a zero volume $\rmd X(t)$.
 \begin{defn}\label{i2} The caustic at time $t$ is defined to be the set,
\[ C_t = \left\{ x: \quad\det\left(\frac{\partial X(t)}{\partial x_0}\right) = 0\right\}. \]
\end{defn}

 Assume that $x$ has $n$ real pre-images,
\[\Phi_t^{-1}\left\{x\right\} =
\left\{x_0(1)(x,t),x_0(2)(x,t),\ldots,x_0(n)(x,t)\right\},\]  where
each $x_0(i)(x,t)\in\mathbb{R}^d$. Then the Feynman-Kac formula and
Laplace's method in infinite dimensions give for a non-degenerate
critical point, \begin{equation}\label{useries}u^{\mu}(x,t)=
\sum\limits_{i=1}^n \theta_i
\exp\left(-\frac{S_0^i(x,t)}{\mu^2}\right),
\end{equation} where
$S_0^i(x,t) :=
S_0\left(x_0(i)(x,t)\right)+A\left(x_0(i)(x,t),x,t\right),$
and $\theta_i$ is an asymptotic series in $\mu^2$. An asymptotic
series in $\mu^2$  can also be found for
$v^{\mu}(x,t)$ \cite{Truman5}.
Note that $\mathcal{S}_t(x) = \min
\{S_0^i(x,t):i=1,2,\ldots,n\}$.

\begin{defn}\label{i3}
The Hamilton-Jacobi level surface is the set,
\[H_t^c = \left\{ x:\quad S_0^i(x,t) =c \mbox{ for
some }i\right\}.\]
The zero level surface $H_t^0$ includes
the wavefront $\mathcal{W}_t$.
\end{defn}
As $\mu\rightarrow 0$, the dominant term in the expansion (\ref{useries}) comes
from the minimising $x_0(i)(x,t)$ which we denote $\tilde{x}_0(x,t)$.  Assuming $\tilde{x}_0(x,t)$ is unique, we obtain the inviscid limit of
the Burgers  fluid velocity
$v^0(x,t) = \dot{\Phi}_t\left(\tilde{x}_0(x,t)\right).$

If the minimising pre-image $\tilde{x}_0(x,t)$ suddenly changes value between two pre-images $x_0(i)(x,t)$ and $x_0(j)(x,t)$, a jump discontinuity will occur in $v^0(x,t)$, the inviscid limit of the Burgers fluid velocity.
There are  two distinct ways in which the minimiser can change; either two pre-images coalesce and disappear (become complex), or the minimiser switches between two pre-images at the same action value. The first of these occurs as $x$ crosses the caustic and when the minimiser disappears the caustic is said to be cool. The second occurs as $x$ crosses the Maxwell set and again, when the minimiser is involved the Maxwell set is said to be cool.

 \begin{defn}\label{m1}
The Maxwell set is,
\begin{eqnarray*}
M_t & = & \left\{x:\, \exists \,x_0,\check{x}_0\in\mathbb{R}^d \mbox{ s.t. }\right.\\
&   &\quad\left.
x=\Phi_t(x_0)=\Phi_t(\check{x}_0), \,x_0\neq \check{x}_0 \mbox{ and }
 \mathcal{A}(x_0,x,t)=\mathcal{A}(\check{x}_0,x,t)
\right\}.
\end{eqnarray*}
 \end{defn}

\begin{eg}[The generic Cusp] Let $V(x,y)=0,$ $k_t(x,y)=0$ and $S_0(x_0,y_0)=x_0^2 y_0/2$.
This initial condition leads to the \emph{generic Cusp}, a semicubical parabolic caustic shown in Figure \ref{cusp}. The caustic $C_t$ (long dash) is given by,
\[
 x_t(x_0)  =  t^2 x_0^3,\quad
 y_t(x_0)  = \hf[3]tx_0^2-\frac{1}{t}.
\]
 The zero level surface $H_t^0$ (solid line) is,
\[
 x_{(t,0)}(x_0)   =  \hf[x_0]\left(1\pm\sqrt{1-t^2 x_0^2}\right),\quad
 y_{(t,0)}(x_0)  =  \frac{1}{2t}\left(t^2x_0^2-1\pm
\sqrt{1-t^2x_0^2}\right)\]
and the Maxwell set $M_t$ (short dash) is $x=0$ for  $y>-1/t.$

\begin{figure}[h!]
\centering
\begin{picture}(5.5,5.5)
    \put(0,0){\resizebox{5.5cm}{!}{\includegraphics{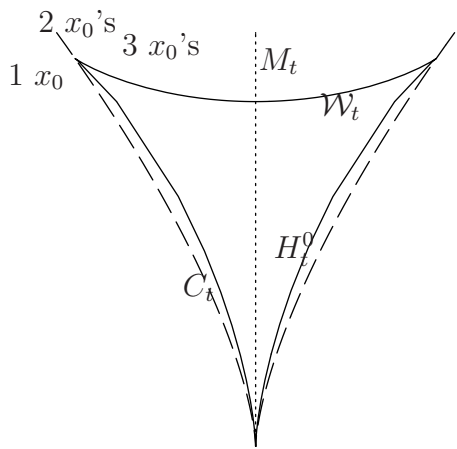}}}
    \put(1.8,2){$C_t$}
    \put(-0.5,4.8){$1$ $x_0$}
    \put(1,5.2){$3$ $x_0$'s}
    \put(-0.1,5.5){$2$ $x_0$'s}
    \put(3.6,4.4){$\mathcal{W}_t$}
    \put(3,2.5){$H_t^0$}
    \put(2.8,5){$M_t$}
  \end{picture}
\caption{Cusp and Tricorn.}\label{cusp}
\end{figure}
\end{eg}

\noindent\emph{Notation:} Throughout this paper
$x,x_0,x_t,x_{(t,c)}$ etc will denote vectors, where normally $x=\Phi_t (x_0)$.
Cartesian coordinates of these will be indicated using a sub/superscript where
relevant; thus $x=(x_1,x_2,\ldots,x_d)$,
$x_0=(x_0^1,x_0^2,\ldots,x_0^d)$ etc. The only exception  will be
 in discussions of explicit examples in two and three dimensions when
we will use $(x,y)$ and $(x_0,y_0)$ etc to denote the vectors.

\section{Some background}

We begin by summarising some results of Davies, Truman and Zhao
(DTZ) \cite{Davies2,Davies3}. Following equation (\ref{action}), let
the stochastic action be defined,
\[
A(x_0,p_0,t) =  \hf\int_0^t \dot{X}(s)^2\rmd s
 -\int_0^t\Bigg[ V(X(s))\rmd s +\epsilon k_s(X(s))\rmd
W_s\Bigg],\]
where
$X(s)=X(s,x_0,p_0)\in\mathbb{R}^d$ and,
\[\rmd\dot{X}(s) = -\nabla V(X(s))\rmd s
-\epsilon \nabla k_s(X(s))\rmd W_s,\quad X(0)=x_0,
\quad\dot{X}(0)=p_0,\] for $ s\in[0,t]$ with
$x_0,p_0\in\mathbb{R}^d$. We assume $X(s)$ is $\mathcal{F}_s$
measurable and unique where $\mathcal{F}_s$ is the sigma algebra
generated by $X(u)$ up to time $s$. It follows from Kunita
\cite{Kunita}:

\begin{lemma} \label{i8}Assume $S_0,V\in
C^2$ and $k_t\in C^{2,0}$, $\nabla V,\nabla k_t$ Lipschitz with Hessians
$\nabla^2 V,\nabla^2 k_t$ and all second derivatives with respect to
space variables of $V$ and $k_t$ bounded. Then for $p_0$, possibly
$x_0$ dependent,
\[\frac{\partial A}{\partial
x_0^{\alpha}}(x_0,p_0,t) =\dot{X}(t)\cdot \frac{\partial
X(t)}{\partial x_0^{\alpha}}-\dot{X}_{\alpha}(0),\qquad\alpha
=1,2,\ldots,d.\]
\end{lemma}

 Methods of Kolokoltsov et al
\cite{Kolokoltsov}  guarantee that for small $t$ the map $p_0\mapsto
X(t,x_0,p_0)$ is onto for all $x_0$.

Therefore, we can define
$A(x_0,x,t)
:=A(x_0,p_0(x_0,x,t),t)$ where
$p_0=p_0(x_0,x,t)$ is the random minimiser (which we assume to be unique) of
$A(x_0,p_0,t)$ when $X(t,x_0,p_0)=x$. Thus, the
stochastic action corresponding to the initial momentum $\nabla
S_0(x_0)$ is $\mathcal{A}(x_0,x,t) := A(x_0,x,t)+S_0(x_0).$

\begin{theorem}\label{sflow} If $\Phi_t$ is the stochastic flow map then,
\[\Phi_t(x_0) = x \quad\Leftrightarrow\quad\frac{\partial}{\partial x_0^{\alpha}}
\left[\mathcal{A}(x_0,x,t)\right]=0,\qquad \alpha=1,2,\ldots,d.\]
\end{theorem}

The Hamilton-Jacobi level surface $H_t^c$ is
found by eliminating $x_0$ between,
\[\mathcal{A}(x_0,x,t) =c, \qquad
\frac{\partial\mathcal{A}}{\partial x_0^{\alpha}}(x_0,x,t) =0 \quad
\alpha=1,2,\ldots, d.\] Alternatively, if we eliminate $x$ to give
an expression in $x_0$, we have the pre-level surface
$\Phi_t^{-1}H_t^c$. Similarly the caustic $C_t$ (and pre-caustic
$\Phi_t^{-1}C_t$) are obtained by eliminating $x_0$ (or $x$) between,
\[\det\left(
\frac{\partial^2\mathcal{A}} {\partial x_0^{\alpha}\partial
x_0^{\beta}}(x_0,x,t) \right)_{\alpha,\beta=1,2,\ldots, d} =0, \qquad
\frac{\partial\mathcal{A}}{\partial
x_0^{\alpha}}(x_0,x,t) =0\quad \alpha = 1,2,\ldots, d.\] These
pre-images are calculated algebraically and are not necessarily the
topological inverse images of the surfaces $C_t$ and $H_t^c$ under
$\Phi_t$.

The caustic surface can be parameterised using its pre-image by
applying the stoch\-astic flow map (a pre-parameterisation). This
allows  us to control the domain of the pre-images and in particular
restrict them to real values. If we can locally solve the equation
of the pre-caustic to give, $x_0^1=\lambda_1,$
$x_0^2=\lambda_2,\ldots,$ $ x_0^{d-1}=\lambda_{d-1},$
$x_0^d=x_{0,\mathrm{C}}^d\left(\lambda_1,\lambda_2,\ldots,\lambda_{d-1}\right),$
then the pre-parameterisation of the caustic is $ x_t(\lambda):=
\Phi_t\left(\lambda,x_{0,\mathrm{C}}^d(\lambda)\right)$ where
$\lambda = (\lambda_1,\ldots,\lambda_{d-1})\in\mathbb{R}^{d-1}$

We next outline a one dimensional analysis first described by
Reynolds, Truman and Williams (RTW) \cite{Truman6}.
\begin{defn}\label{i16}
    The $d$-dimensional  flow map $\Phi_{t}$ is globally reducible if for
    any
    $x=(x_{1},x_{2},\ldots ,x_{d})$ and $x_{0}=(x_{0}^{1},x_{0}^{2},\ldots,
    x_{0}^{d})$ where $x=\Phi_{t}(x_{0}),$ it is possible to write
    each coordinate $x_0^\alpha$ as a function of the lower
    coordinates. That is,
    \begin{equation}\label{ie16}
    x=\Phi_{t}(x_{0}) \quad \Rightarrow\quad
        x_{0}^{\alpha}=x_{0}^{\alpha}(x,x_{0}^{1},x_{0}^{2},\ldots, x_{0}^{\alpha-1},t)\mbox{ for $ \alpha=d,d-1,\ldots ,2.$}
    \end{equation}

\end{defn}
Therefore, using Theorem \ref{sflow}, the flow map is globally
reducible if we can find a chain of $C^{2}$ functions
$x_{0}^{d},x_{0}^{d-1},\ldots,x_0^2$ such that,
\begin{eqnarray*}
     x_{0}^{d}=x_{0}^{d}(x,x_{0}^{1},x_{0}^{2},\ldots
   x_{0}^{d-1},t)
&\Leftrightarrow &
\frac{\partial\mathcal{A}}{\partial x_{0}^{d}}(x_{0},x,t)=0,\\
 x_{0}^{d-1}=x_{0}^{d-1}(x,x_{0}^{1},x_{0}^{2},\ldots
x_{0}^{d-2},t) &\Leftrightarrow &
\frac{\partial\mathcal{A}}{\partial
x_{0}^{d-1}}(x_{0}^{1},x_{0}^{2},\ldots,x_{0}^{d}(\ldots),x,t)=0,\\
& \vdots &\\
 x_{0}^{2}=x_{0}^{2}(x,x_{0}^{1},t) &\Leftrightarrow &
 \end{eqnarray*}
\[\qquad\qquad\qquad\qquad
\qquad \qquad  \,\,\,\,\frac{\partial\mathcal{A}}{\partial
x_{0}^{2}}(x_{0}^{1},x_{0}^{2},x_{0}^{3}(x,x_{0}^{1},x_{0}^{2},t),\ldots,x_{0}^{d}(\ldots),x,t)=0.
\]
This requires that no roots are repeated to ensure that none of the
second derivatives of $\mathcal{A}$ vanish. We assume also that
there is a favoured ordering of coordinates and a corresponding
decomposition of $\Phi_{t}$ which allows the non-uniqueness to be
reduced to the level of the $x_{0}^{1}$ coordinate. This assumption
appears to be quite restrictive. However, local reducibility at $x$
follows from the implicit function theorem and some mild assumptions
on the derivatives of $\mathcal{A}$ (see \cite{Neate}).

\begin{defn}\label{i17}
If $\Phi_{t}$ is globally reducible then  the reduced
    action function is the univariate function obtained by evaluating the
     action with equations (\ref{ie16}),
    \[f_{(x,t)}(x_0^1):=f(x_{0}^{1},x,t) =\mathcal{A}
    (x_{0}^{1},x_{0}^{2}(x,x_{0}^{1},t),x_{0}^{3}(\ldots),\ldots,x,t).\]
\end{defn}

\begin{lemma} \label{i18}If $\Phi_t$ is globally reducible, modulo the above assumptions,
\begin{eqnarray*}\lefteqn{\left|
    \det
    \left.\left(\frac{\partial^{2}\mathcal{A}}{(\partial
    x_{0})^{2}}({x}_{0},{x},t)\right)
    \right|_{x_0=(x_0^1,x_0^2(x,x_0^1,t),\ldots,x_0^d(\ldots))}\right|}\\
    &  = & \prod\limits_{\alpha=1}^{d}\left|\left[
    \left(\frac{\partial}{\partial
 x_{0}^{\alpha}}\right)^2\!\!\! \mathcal{A}
 (x_{0}^{1},\ldots,x_{0}^{\alpha},x_0^{\alpha+1}(\ldots),\ldots
    ,x_{0}^{d}(\ldots),{x},t)\right]_{\begin{array}{c}
   { \scriptstyle x_0^2=x_0^2(x,x_0^1,t)}\\[-1ex]
     {\scriptstyle{\vdots}}\\[-1ex]
   {\scriptstyle x_0^{\alpha} = x_0^{\alpha}(\ldots)}
    \end{array}}\right|
\end{eqnarray*}
   where the first term is $f_{({x},t)}''(x_{0}^{1})$ and the last
    $d-1$ terms are nonzero.
\end{lemma}

\begin{theorem}\label{i19}
Let the classical mechanical flow map $\Phi_t$ be globally
reducible.   Then:
    \begin{enumerate}
    \item $ f_{(x,t)}' (x_{0}^{1})=0$
     and the  equations (\ref{ie16})
     \emph{$\Leftrightarrow x=\Phi_{t} (x_{0}),$}
    \item
    $ f_{(x,t)}' (x_{0}^{1})
    = f_{(x,t)}'' (x_{0}^{1})=0$
    and the  equations (\ref{ie16})\\
    \emph{$\Leftrightarrow x=\Phi_{t}(x_{0})$} is such that the number of real
    solutions $x_{0}$ changes.
    \end{enumerate}
\end{theorem}
\begin{cor}\label{c2}
Let $x_{t}(\lambda)$ denote the pre-parameterisation of the caustic with
$\lambda
    =(\lambda_{1},\lambda_{2},\ldots,\lambda_{d-1})\in\mathbb{R}^{d-1}$. Then
$f'_{(x_{t}(\lambda),t)}(\lambda_{1}) =
    f''_{(x_{t}(\lambda),t)}(\lambda_{1})=0.$
\end{cor}

\begin{figure}[ht]
\setlength{\unitlength}{8mm}
\noindent\begin{tabular}{ccc} \textbf{Before Caustic}& \textbf{On
Cool Caustic} & \textbf{Beyond Caustic}\\
\begin{picture}(4.5,3) \put(0,0){\resizebox{36mm}{!}{\includegraphics{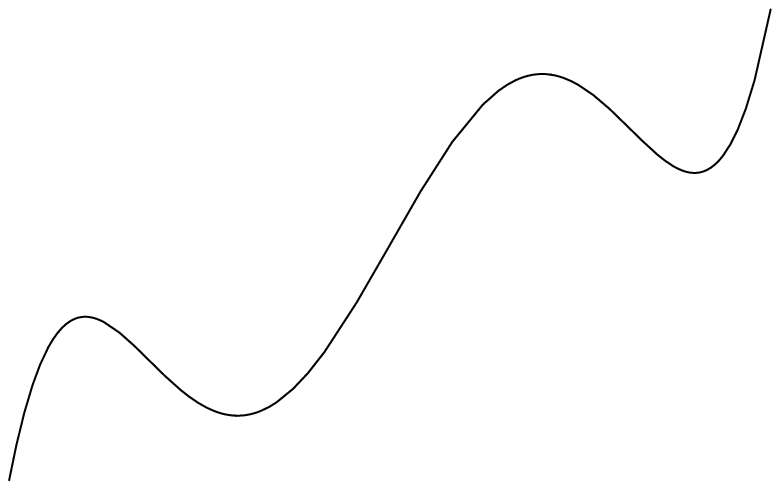}}}
\put(0.2,1.2){\scriptsize$x_0^1(1)$}
\put(0.4,0.15){\scriptsize$x_0^1(2)=\tilde{x}^1_0(x,t)$} \put(2.65,2.5){\scriptsize$x_0^1(3)$}
\put(3.6,1.4){\scriptsize $x_0^1(4)$}
\put(1.2,0.65){\scriptsize $(a)$}
\put(3.7,2){\scriptsize $(b)$}
\end{picture}&
\begin{picture}(4.5,3) \put(0,0){\resizebox{36mm}{!}{\includegraphics{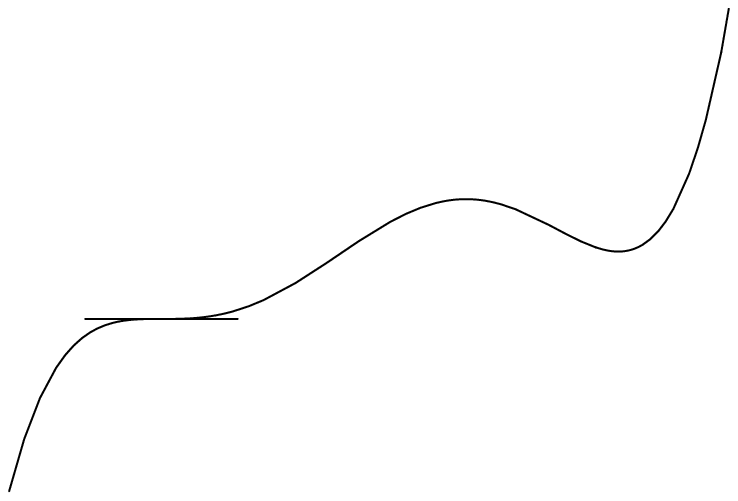}}}
\put(0.8,1.2){\scriptsize $(a)$}
\end{picture}&
\begin{picture}(4.5,3)
\put(0,0){\resizebox{36mm}{!}{\includegraphics{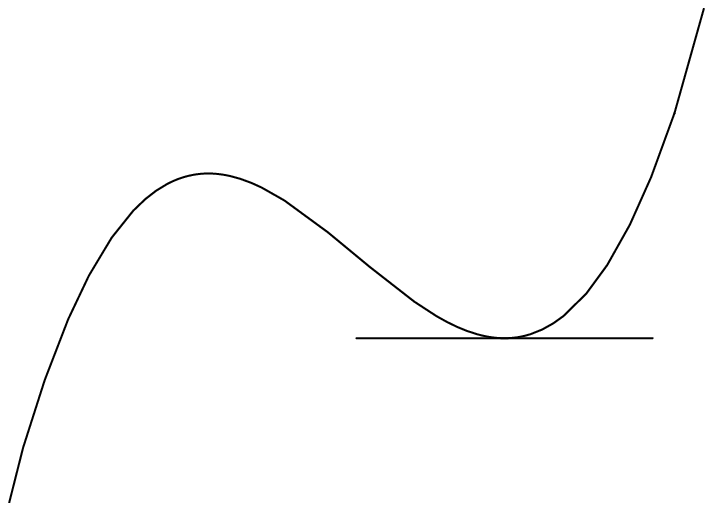}}}
\put(2.6,1.1){\scriptsize $(b)$}
\end{picture}\\
\emph{ Minimiser at }   & \emph{ Two $x_0^1$'s coalesce to} &
\emph{ Minimiser jumps.} \\
      $x_0^1(2)(x,t)=\tilde{x}^1_0(x,t)$.              &\emph { form  a point of inflexion.} &
\end{tabular}
 \caption{The graph of $f_{(x,t)}(x_0^1)$ as $x$ crosses the caustic.}\label{cf1}
\end{figure}

 Consider an example where
for $x$ on one side of the caustic there are four real critical
points on $f_{(x,t)}(x_0^1)$  enumerated $x_0^1(i)(x,t)$ for $i=1$
to $4$, and denote the minimising critical point by
$\tilde{x}^1_0(x,t)$. Figure \ref{cf1} illustrates how the minimiser
jumps from $(a)$ to $(b)$ as $x$ crosses the caustic if the point of
inflexion at $x_0^1=\lambda$ is the global minimiser. In this case
the caustic at $x_t(\lambda)$ is said to be cool.

\begin{figure}[h!]
 \setlength{\unitlength}{0.85cm}
\centering
\begin{tabular}{ccc} \textbf{Before $M_t$}& \textbf{On
Cool $M_t$} & \textbf{Beyond $M_t$}\\
\begin{picture}(4.5,3) \put(0,0){\resizebox{40.5mm}{!}{\includegraphics{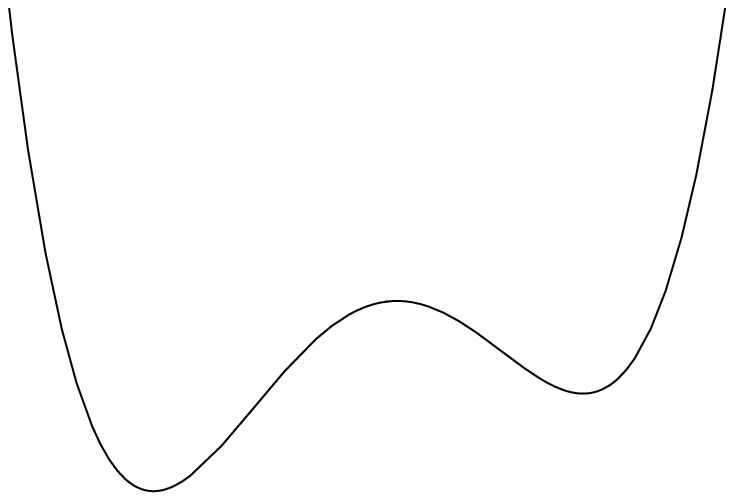}}}
 \put(0.9,0.4){\small{$x_0^1$}}
 \put(3.4,0.9){\small{$\check{x}_0^1$}}\end{picture}&
\begin{picture}(4.5,3) \put(0,0){\resizebox{40.5mm}{!}{\includegraphics{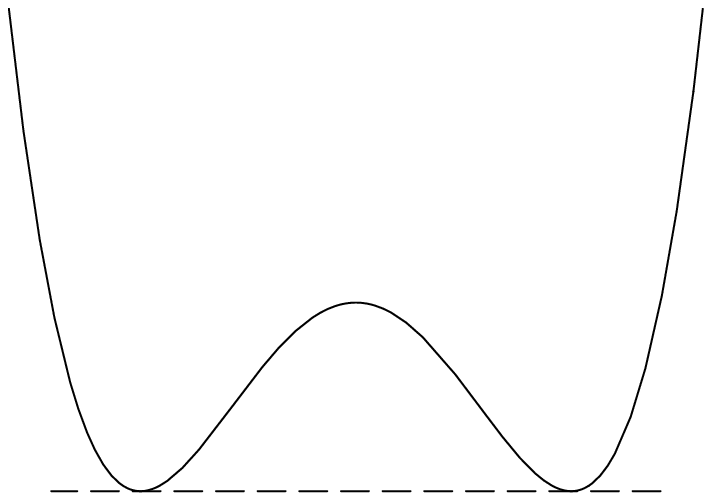}}}
  \put(0.9,0.4){\small{$x_0^1$}}
 \put(3.4,0.4){\small{$\check{x}_0^1$}}\end{picture}&
\begin{picture}(4.5,3)
\put(0,0){\resizebox{40.5mm}{!}{\includegraphics{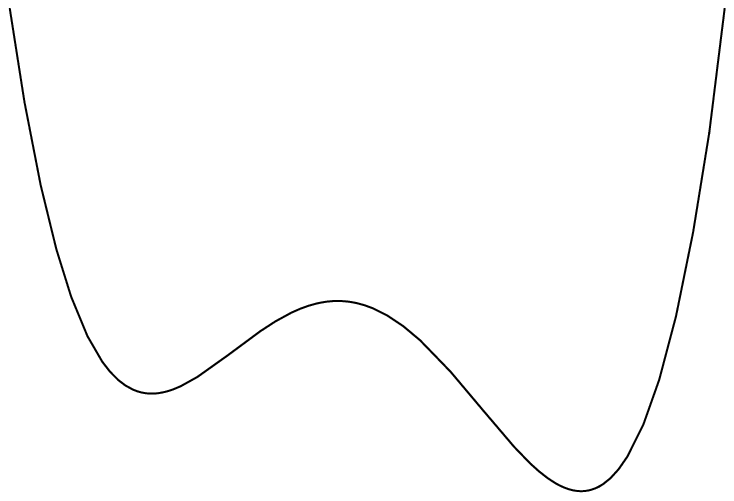}}}
  \put(1,0.9){\small{$x_0^1$}}
 \put(3.5,0.4){\small{$\check{x}_0^1$}}
 \end{picture}\\
\emph{ Minimiser at $x_0^1$.}   & \emph{ Two $x_0$'s at same level.}
& \emph{ Minimiser jumps.}
\end{tabular}
\caption{The graph of $f_{(x,t)}(x_0^1)$ as $x$ crosses the Maxwell
set.}\label{mf1}
\end{figure}

We can also consider the Maxwell set in terms of the reduced action
function. The Maxwell set corresponds to values of $x$ for which
$f_{(x,t)}(x_0^1)$ has two critical points at the same action value.
If each of this pair of critical points also minimises the reduced
action, then the inviscid limit of the solution to the Burgers
equation will jump as shown in Figure \ref{mf1} and the Maxwell set
will be described as cool. Note that a Maxwell set can only exist in
a region with three or more real pre-images if the reduced action is
continuous.

We now use the reduced action function to find the Maxwell set.
Instead of finding the Maxwell set directly it is easier to find the
Maxwell-Klein set where Definition \ref{m1} is changed to allow
$x_0,\check{x}_0\in\mathbb{C}^d$ \cite{Neate}.

\begin{theorem}\label{m10}
Let the reduced action function $f_{(x,t)}(x_0^1)$ be polynomial
in all space variables. Then the set of possible discontinuities
for a $d$-dimensional Burgers fluid velocity field in the inviscid
limit is the double discriminant,
 $$D(t):=D_c\left\{D_{\lambda}\left(f_{(x,t)}(\lambda)-c\right)\right\}=0,$$
 where $D_x(p(x))$ is the discriminant of the polynomial $p$ with
 respect to $x$.
\end{theorem}

\begin{theorem}\label{m13}
The double discriminant $D(t)$ factorises as,
\[D(t) = b_0^{2m-2}\cdot\left(C_t\right)^3\cdot \left(B_t\right)^2\]
where $B_t=0$ is the equation of the Maxwell-Klein set, $C_t=0$
is the equation of the caustic and $b_0$ is some function of $t$ only; both $B_t$ and $C_t$
 are algebraic in
$x$ and $t$.
\end{theorem}
Theorem \ref{m13} gives the Maxwell-Klein set as an algebraic
surface. It is then necessary to extract the Maxwell set from this
by establishing which points have real pre-images. Alternatively, we
can find the Maxwell set via a pre-parameterisation which allows us
to restrict the pre-images to be real; for this we need the
pre-Maxwell set which can also be found using the reduced action
function.
\begin{theorem}\label{m17}
The pre-Maxwell set is given by the discriminant,
\[D_{\check{x}_0^1}\left(\frac{f_{(\Phi_t (x_0),t)}(x_0^1)-f_{(\Phi_t (x_0),t)}(\check{x}_0^1)}{(x_0^1-\check{x}_0^1)^2}\right)=0.\]
\end{theorem}

\section{Geometric results}
We now summarise geometric results of DTZ  and also \cite{Neate3}.
Assume that $A(x_0,x,t)$ is $C^4$ in space variables with
$\det\left( \frac{\partial^2\mathcal{A}}{\partial
x_0^{\alpha}\partial x^{\beta}}\right)\neq 0.$

\begin{defn}
A curve $x=x(\gamma)$, $\gamma \in N(\gamma_0,\delta)$, is said to
have a generalised cusp at $\gamma= \gamma_0$, $\gamma$ being an
intrinsic variable such as arc length, if
$x'( \gamma_0)=0.$
\end{defn}
\begin{lemma}\label{s11}
Let $\Phi_t$ denote the stochastic flow map and
$\Phi_t^{-1}\Gamma_t$ and $\Gamma_t$ be some surfaces where if
$x_0\in\Phi_t^{-1}\Gamma_t$ then $x=\Phi_t (x_0)\in\Gamma_t$. Then,
$\Phi_t$ is a differentiable map from $\Phi_t^{-1}\Gamma_t$ to
$\Gamma_t$ with Frechet derivative,
\[(D\Phi_t)(x_0) = \left(-\frac{\partial^2\mathcal{A}}{\partial x\partial x_0}(x_0,x,t)\right)^{-1}
\left(\frac{\partial^2\mathcal{A}}{(\partial x_0)^2}(x_0,x,t)\right).
\]
\end{lemma}
\begin{lemma}[2 dims]\label{s12}
Let $x_0(s)$ be any two dimensional intrinsically parameterised
curve, and define
$x(s)=\Phi_t(x_0(s)).$
Let $e_0$ denote the zero eigenvector of $(\partial^2\mathcal{A}/(\partial x_0)^2)$ and
assume that $\ker(\partial^2\mathcal{A}/(\partial x_0)^2)=\langle e_0\rangle$. Then, there is a generalised
cusp on $x(s)$ when $s=\sigma$ if and only if either:
\begin{enumerate}
\item there is a generalised cusp on $x_0(s)$ when $s=\sigma$; or,
\item $x_0(\sigma)\in\Phi_t^{-1}C_t$ and the tangent
$\frac{\di x_0}{\di s}(s)$ at $s=\sigma$ is parallel to $e_0$.
\end{enumerate}
\end{lemma}

\begin{prop}The normal to the pre-level
surface is,
\[n_{\mathrm{H}}(x_0) =
-\left(\frac{\partial^2\mathcal{A}}{(\partial x_0)^2}\right)
\left(\frac{\partial^2\mathcal{A}}{\partial x_0\partial x}\right)^{-1}
\dot{X}\left(t,x_0,\nabla S_0(x_0)\right).\]

\end{prop}

\begin{prop}[2 dims]\label{i11}
Assume that at $x_0\in\Phi_t^{-1} H_t^c$ the normal to
$\Phi_t^{-1}H_t^c$ is non-zero, so that the pre-level surface does
not have a generalised cusp at $x_0$. Then the level surface can
only have a generalised cusp at $\Phi_t(x_0)$ if $\Phi_t(x_0)\in
C_t$. Moreover, if $x=\Phi_t(x_0)\in\Phi_t\left\{\Phi_t^{-1}C_t\cap
\Phi_t^{-1}H_t^c\right\}$, the level surface will have a generalised
cusp.
\end{prop}

\begin{eg}[The generic Cusp]
Figure \ref{if2} shows how a point lying on three level surfaces  has three distinct real
pre-images each on a separate pre-level surface. A cusp only occurs on the corresponding level surface when the pre-level surface
intersects the pre-caustic. Thus, provided the normal to the pre-level surface is well defined, a level surface can only have a
cusp on the caustic, but it does not have to be cusped when it meets the caustic.
\begin{figure}[h!]
\centering
\begin{tabular}{cc}
 \begin{picture}(5.5,5.5)
 \put(0,0){\resizebox{55mm}{!}{\includegraphics{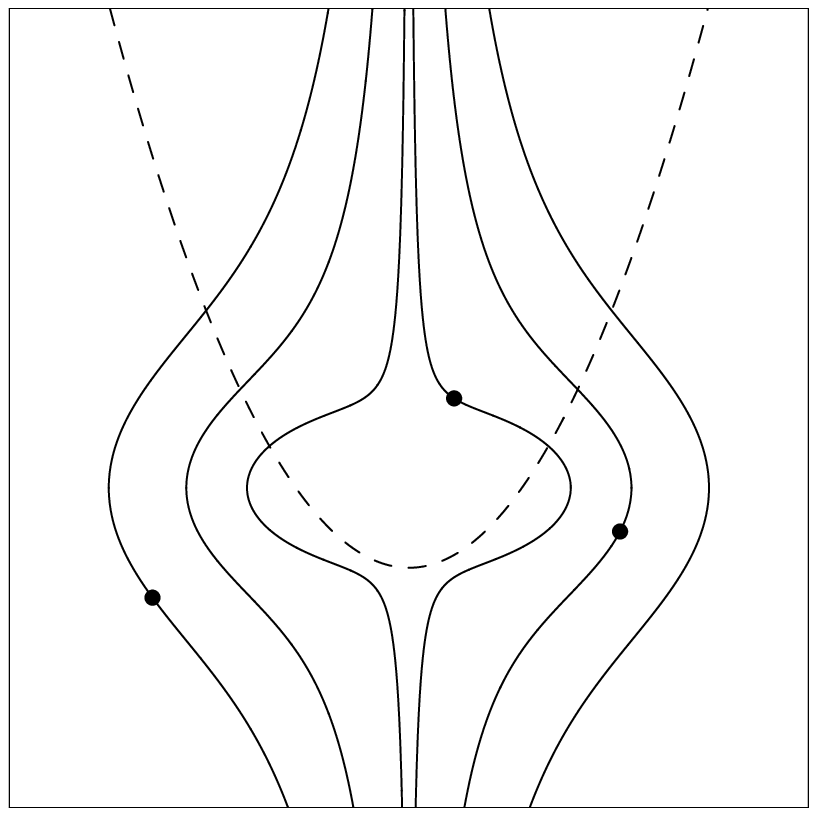}}}
 \put(0.1,0.3){(a)}
 \end{picture}&
\begin{picture}(5.5,5.5)
 \put(0,0){\resizebox{55mm}{!}{\includegraphics{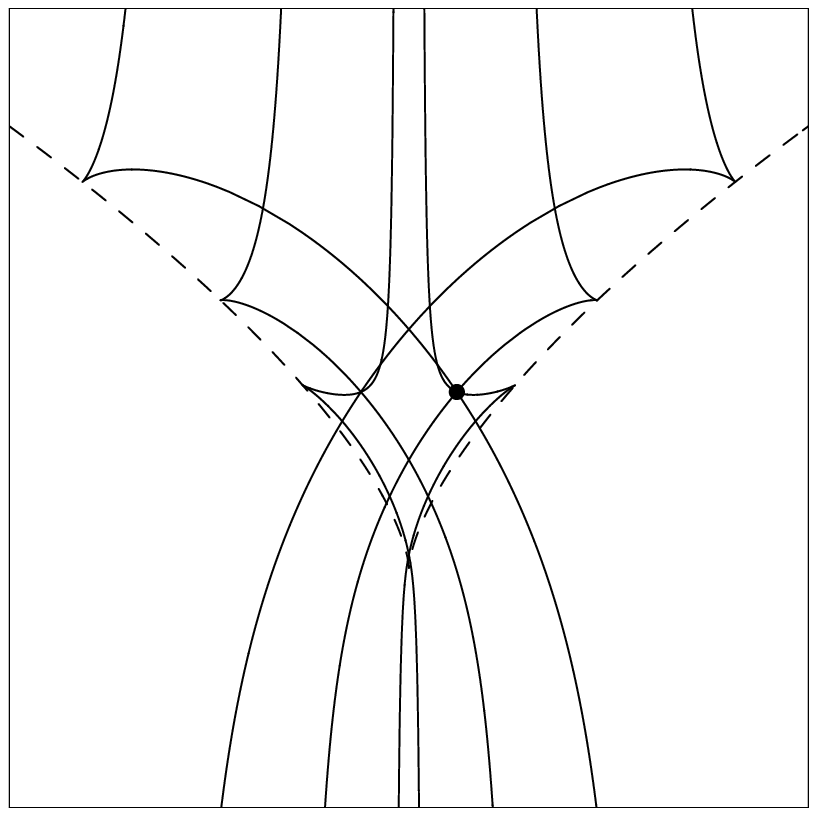}}}
 \put(0.1,0.3){(b)}
 \end{picture}
\end{tabular}
 \caption{(a) The pre-level surface (solid line) and pre-caustic
 (dashed),
 (b) the level surface (solid line) and caustic (dashed), both for the generic Cusp.}
 \label{if2}
 \end{figure}
\end{eg}
\begin{prop}
Assume that $x\in M_t$ corresponds to exactly
two pre-images on the pre-Maxwell set, $x_0$ and $\check{x}_0$. Then the normal to the pre-Maxwell set at $x_0$ is given by,
\begin{eqnarray*}
n_{\mathrm{M}}(x_0) &  = &  - \left(\frac{\partial^2\mathcal{A}}{(\partial
x_0)^2}(x_0,x,t)\right)
\left(\frac{\partial^2\mathcal{A}}{\partial x_0\partial x}(x_0,x,t)\right)^{-1}\cdot\\
& &\qquad
\left(\dot{X}(t,x_0,\nabla S_0(x_0))-\dot{X}(t,\check{x}_0,\nabla S_0(\check{x}_0))\right)_{\displaystyle .}\end{eqnarray*}
\end{prop}
\begin{prop}[2 dims]\label{m24}
Assume that at $x_0\in\Phi_t^{-1}M_t$,
$n_{\mathrm{M}}(x_0)\neq 0$ so that the pre-Maxwell set does not have a generalised cusp. Then, the Maxwell set can only have a cusp at $\Phi_t(x_0)$ if $\Phi_t(x_0)\in C_t$. Moreover, if
$x=\Phi_t(x_0)\in\Phi_t\left\{ \Phi_t^{-1}C_t\cap
\Phi_t^{-1}M_t\right\},$ the Maxwell set will have a generalised
cusp at $x$.
\end{prop}
\begin{cor}[2 dims]\label{m25}
Assuming that $n_{\mathrm{H}}(x_0)\neq 0$ and $n_{\mathrm{M}}(x_0)\neq 0$, then when the pre-Maxwell set intersects the pre-caustic, it touches a pre-level surface.
Moreover, there is a cusp on the pre-Maxwell set which also intersects the same pre-level surface.
\end{cor}

\begin{eg}[The polynomial swallowtail] Let $V=0$, $k_t=0$ and $S_0(x_0,y_0)=x_0^5+x_0^2 y_0$.
From Proposition \ref{m24}, the cusps on the Maxwell set correspond
to the intersections of the pre-curves (points 3 and 6 on Figure \ref{mf4}). But from
Corollary \ref{m25}, the cusps on the Maxwell set also correspond to
the cusps on the pre-Maxwell set (points 2 and 5). Each cusp on the
pre-Maxwell set lies on the same level surface as a point of
intersection between the pre-caustic and pre-Maxwell set
(ie. 3 and 5, 2 and 6).
\begin{figure}[h]
\setlength{\unitlength}{0.76cm}
\centering \begin{tabular}{cc}
 \begin{picture}(7,7)
 \put(2.3,4.75){1}
  \put(0.4,5.3){2}
 \put(1.6,5.9){3}
  \put(4.8,5.3){4}
  \put(6.7,0.1){5}
 \put(5.15,2.95){6}
  \put(0.1,0.3){(a)}
\put(0,0){\resizebox{55mm}{!}{\includegraphics{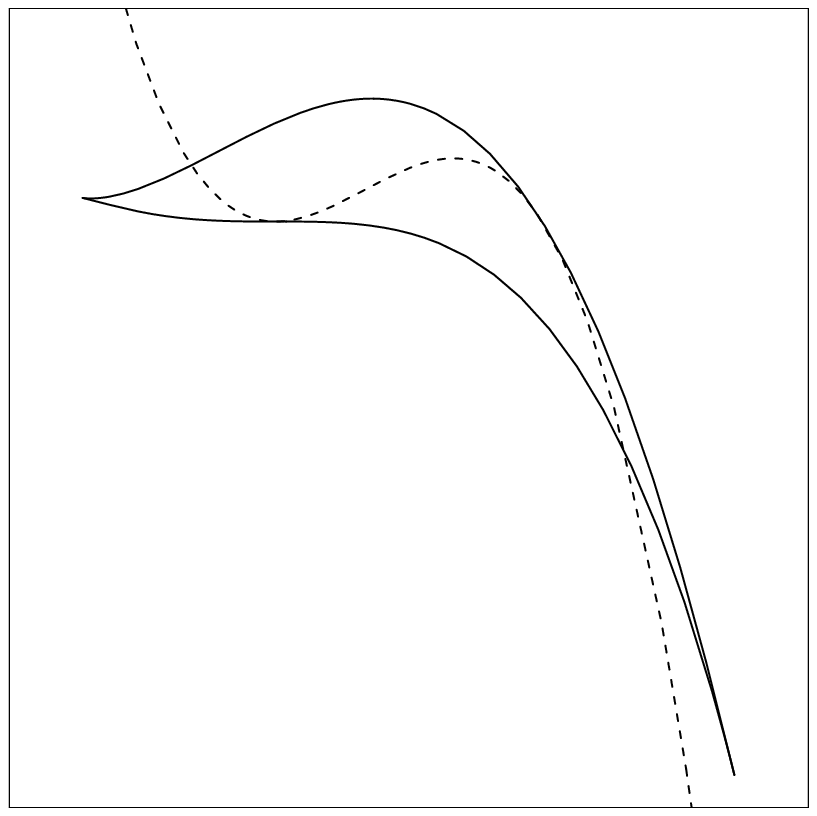}}}
\end{picture}
&
 \begin{picture}(7,7)
 \put(2.7,0.2){1}
  \put(3.2,4.7){2}
 \put(2,2.5){3}
  \put(6.9,6.5){4}
  \put(2.2,2.){5}
 \put(3.6,4.8){6}
  \put(0.1,0.3){(b)}
\put(0,0){\resizebox{55mm}{!}{\includegraphics{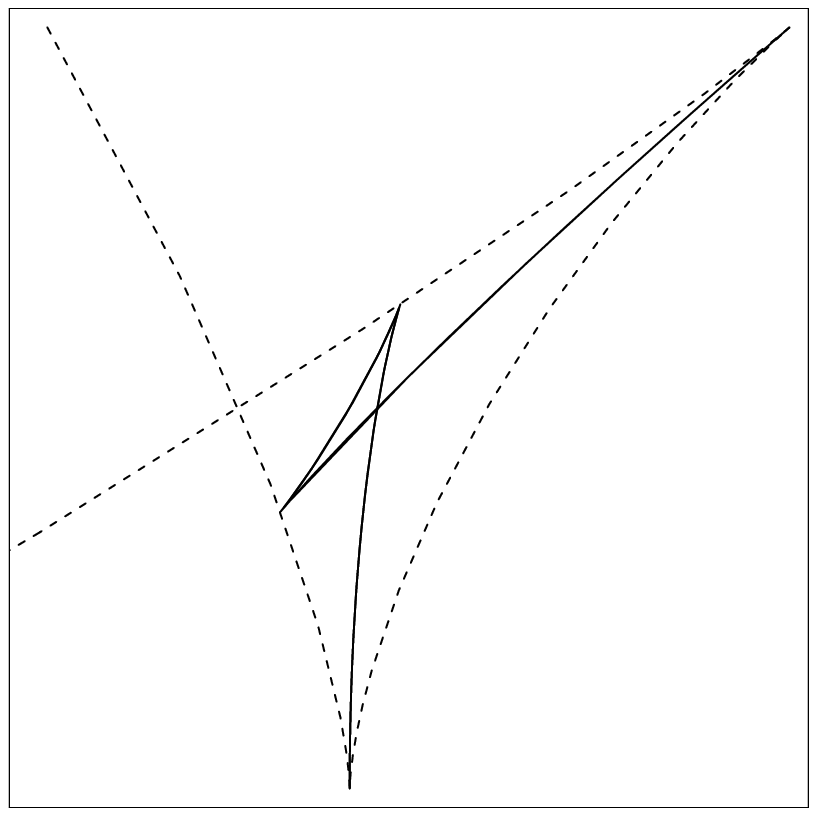}}}
\end{picture}
\end{tabular}
\caption{(a) The pre-Maxwell set (solid line) and pre-caustic
 (dashed),
 (b) the Maxwell set (solid) and caustic (dashed), for the polynomial swallowtail.} \label{mf4}
\end{figure}

The Maxwell set terminates when it reaches the cusps on the caustic
(points 1 and 4). These points satisfy the condition for a
generalised cusp but, instead of appearing cusped, the curve stops
and the parameterisation begins again in the sense that it maps back
exactly onto itself. This follows because every point on the Maxwell
set has at least two real pre-images, and so by pre-parameterising
the Maxwell set, we effectively sweep it out twice. All of the
pre-surfaces touch at the cusps on the caustic.

\end{eg}
These results can be extended to three dimensions.

\begin{theorem}[$3$ dims] Let
\[
x \in  \mbox{Cusp}\left(H_t^c\right)
=
\left\{x\in\Phi_t\left(\Phi_t^{-1}C_t\cap\Phi_t^{-1}H^c_t\right),x=\Phi_t(x_0),
n_{\mathrm{H}}(x_0)\neq 0\right\}.\]Then in three dimensions, with probability one, $T_{\mathrm{H}}(x)$ the tangent space to
the level surface at $x$ is at most one dimensional.
\end{theorem}
\begin{theorem}[3 dims]\label{m29}
Let,
\[x\in\mathrm{Cusp}(M_t) =
\left\{
x\in \Phi_t\left(\Phi_t^{-1}C_t\cap\Phi_t^{-1} M_t\right),
 x=\Phi_t(x_0), n_{\mathrm{M}}(x_0)\neq 0
\right\}.\]
 Then in three dimensions, with probability one, $T_{\mathrm{M}}(x)$, the tangent space to the Maxwell set at $x$ is at most one dimensional.
\end{theorem}

\section{Swallowtail perestroikas}
The geometry of a caustic or wavefront
can suddenly change with singularities appearing and disappearing \cite{Arnold3}. We consider the formation or collapse of a swallowtail using some earlier works of Cayley and Klein. Here we provide a summary of results from \cite{Neate}.

\begin{figure}[h]
  \setlength{\unitlength}{1 cm} \centering
 \resizebox{110mm}{!}{ \includegraphics{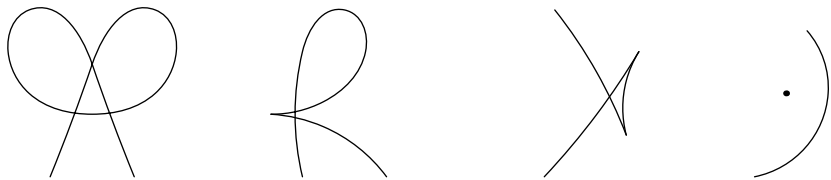}}
  \caption{Cayley's triple points.}
  \label{pic3}
  \end{figure}

In Cayley's work on plane algebraic curves, he describes the possible
triple points of a curve \cite{Salmon} by considering the
collapse of systems of double points which would lead to the
existence of three tangents at a point. The four possibilities are
shown in Figure \ref{pic3}. The systems will collapse to form a
triple point with respectively,
three real distinct tangents,
 three real tangents with two coincident,
 three real tangents all of which are coincident,
or one real tangent and two complex tangents.
We are interested in the interchange
between the last two cases which Felix Klein investigated \cite{Hwa,Klein}.

As indicated in Section 2,  we often parameterise the caustic and
level surfaces using a pre-parameterisation in which we restrict the
parameter to be real to only consider points with real pre-images.
This does not allow there to be any isolated double points on these
curves. We now let the pre-parameter vary throughout the complex
plane and consider when this maps to real points. We begin  with a
family of curves  of the form ${x}_t(\lambda)
=(x_{t}^1(\lambda),x_t^2(\lambda))$ where each
$x_t^{\alpha}(\lambda)$ is real analytic in $\lambda$. If
Im$\left\{{x}_t(a+\rmi\eta)\right\}=0$,
 it follows that ${x}_t(a+\rmi\eta)={x}_t(a-\rmi\eta),$
  so this is a
  ``\emph{complex double point}'' of the curve ${x}_t(\lambda)$.

 \begin{prop}\label{swall}
      If a swallowtail on the curve ${x}_{t}(\lambda)$
      collapses to a point
      where $\lambda=\tilde{\lambda}$ when $t=\tilde{t}$ then
      $\frac{\rmd {x}_{\tilde{t}}}{\rmd \lambda}(\tilde{\lambda})=
      \frac{\rmd^{2} {x}_{\tilde{t}}}{\rmd \lambda^{2}}(\tilde{\lambda})=0.$
  \end{prop}

\begin{prop}\label{comp}
Assume that there exists a neighbourhood of $\tilde{\lambda}\in\mathbb{R}$
such that $\frac{\rmd {x}_{t}^{\alpha}}{\rmd \lambda}(\lambda)\neq 0$
for $t\in(\tilde{t}-\delta,\tilde{t})$ where $\delta>0$.
    If a complex double point joins the curve ${x}_{t}(\lambda)$
    at $\lambda=\tilde{\lambda}$ when $t=\tilde{t}$ then
    $\frac{\rmd {x}_{\tilde{t}}}{\rmd \lambda}(\tilde{\lambda})=
      \frac{\rmd^{2} {x}_{\tilde{t}}}{\rmd \lambda^{2}}(\tilde{\lambda})=0.$

\end{prop}

    These give a necessary condition for the
formation or destruction of a swallowtail,  and for complex double
points to join or leave the main curve.
 \begin{defn}
 A  family of parameterised curves $x_t(\lambda)$, (where $\lambda$ is some intrinsic parameter) for which
 $\frac{\rmd x_{\tilde{t}}}{\rmd \lambda}(\tilde{\lambda})=
 \frac{\rmd^2 x_{\tilde{t}}}{\rmd \lambda^2}(\tilde{\lambda})=0$
  is said to have  a point of swallowtail perestroika when $\lambda=\tilde{\lambda}$ and $t=\tilde{t}$.
 \end{defn}
As with generalised cusps,
we have not ruled out further degeneracy at these points.
Moreover, as Cayley highlighted, these points are not cusped and are barely distinguishable from an ordinary point of the curve \cite{Salmon}.

We now apply these ideas to the caustic where $x_t(\lambda)$ will
denote the pre-parameterisation. The ``\emph{complex caustic}" is
found by allowing the parameter $\lambda$  to vary over the complex
plane. We are interested in the complex double points if they join
the main caustic at some finite critical time $\tilde{t}$ where
$\eta_{t}\rightarrow 0$ as $t\uparrow \tilde{t}$; at such a point a
swallowtail can develop.

 \begin{eg}Let $V(x,y)= 0, k_t(x,y)\equiv 0$ and
    $S_{0}(x_{0},y_0) = x_{0}^{5}+x_{0}^{6}y_{0}.$
The caustic  has no cusps for times $t<\tilde{t}$ and two cusps for times
    $t>\tilde{t}$ where,
    $\tilde{t} = \frac{4}{7}\sqrt{2} \left(\frac{33}{7}\right)^{\frac{3}{4}}
    =2.5854\ldots$

At the critical time $\tilde{t}$ the caustic has a point of swallowtail perestroika
as shown in Figure  \ref{dash2}.  There are   five complex double points before the critical time and four
afterwards. The remaining complex double points do not join the main
caustic and so do not influence its behaviour for real times.
\begin{figure}[h!]
 \centering\resizebox{120mm}{!}{\includegraphics{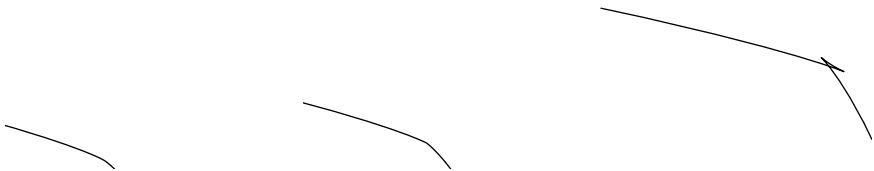}}
\caption{Caustic plotted at corresponding times.}
\label{dash2}
\end{figure}
\end{eg}

Unsurprisingly, these phenomena are not restricted to caustics. There is an interplay  between the level surfaces
and the caustics, characterised by their pre-images.
\begin{prop}\label{s13}
    Assume that in two dimensions at
    $x_0\in\Phi_t^{-1}H_t^c\cap\Phi_t^{-1}C_t$
    the normal to the pre-level surface $n_{\mathrm{H}}(x_0)\neq 0$ and the normal
    to the pre-caustic $n_{\mathrm{C}}(x_0)\neq 0$ so that the pre-caustic
    is not cusped at $x_0$.
    Then  $n_{\mathrm{C}}(x_0)$ is parallel to $n_{\mathrm{H}}(x_0)$ if and only if
    there is a generalised cusp on the caustic.
    \end{prop}
    \begin{cor}\label{s14}
    Assume that in two dimensions at
    $x_0\in\Phi_t^{-1}H_t^c\cap\Phi_t^{-1}C_t$
    the normal to the pre-level surface $n_{\mathrm{H}}(x_0)\neq 0$.
    Then at $\Phi_t(x_0)$ there
    is a point of swallowtail perestroika on the level surface
    $H_t^c$ if and only if  there is a generalised cusp on the caustic $C_t$ at $\Phi_t(x_0)$.
\end{cor}
 \begin{eg}Let $V(x,y)= 0$, $k_t(x,y)=0$, and
 $S_0(x_0,y_0)=x_0^5+x_0^6y_0.$
Consider the behaviour of the level surfaces through a
point inside the caustic swallowtail at a fixed time as the point is moved through a cusp on the
caustic. This is illustrated in Figure  \ref{cuspls}. Part (a) shows
all five of the level surfaces through the point demonstrating how three
swallowtail level surfaces collapse together at the cusp to form a
single level surface with a point of swallowtail perestroika. Parts
(b) and (c) show how one of these swallowtails collapses on its own
and how its pre-image behaves.
\begin{figure}[h!]
\setlength{\unitlength}{10mm}
\begin{picture}(12.5,5.5)
\put(0.5,0){\resizebox{120mm}{!}{\includegraphics{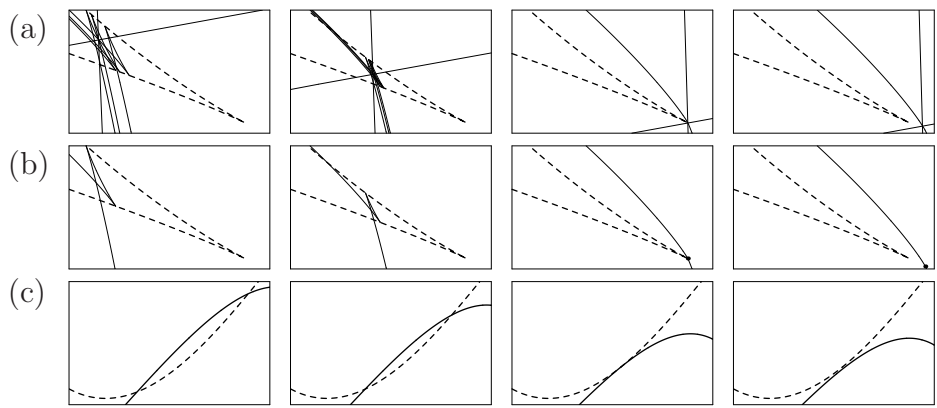}}}
\put(0,1.5){(c)}
\put(0,3.2){(b)}
\put(0,5){(a)}
\end{picture}
\caption{(a) All level surfaces (solid line) through a point as it
crosses the caustic (dashed line) at a cusp, (b) one of these level
surfaces with its complex double point, and (c) its real pre-image.}
\label{cuspls}
 \end{figure}
\end{eg}

\section{Real turbulence}
The geometric results of Section 3 showed that cusps (or in three
dimensions curves of cusps) on the level surfaces occur where the
pre-level surface intersects the pre-caustic. As time passes,  the
cusps or curves of cusps will appear and disappear on the level
surfaces as the pre-curves move.

\begin{defn}\label{t1}
    Real turbulent times are defined to be times $t$ at which there
    exist points where the
    pre-level surface $\Phi_t^{-1}H_t^c$ and pre-caustic $\Phi_t^{-1} C_t$ touch.
\end{defn}
Real turbulent times correspond to times at which there is a change
in the number of cusps or cusped curves on the level surface
$H_{t}^{c}$.

In $d$-dimensions, assuming $\Phi_t$ is
globally reducible, let $f_{({x},t)}(x_{0}^1)$ denote the reduced
action function and $x_t(\lambda)$ denote the pre-parameteris\-ation of the
caustic.

\begin{theorem}\label{t2}
The real turbulent times $t$ are given by the zeros of the zeta
process $\zeta_t^c$ where,
\[\zeta_t^c:=f_{(x_{t}(\lambda),t)}(\lambda_1)-c,\]
$\lambda$ satisfies,
\begin{equation}\label{te1b}\frac{\partial}{\partial \lambda_{\alpha}}f_{(x_{t}(\lambda),t)}(\lambda_1)
=0\quad\mbox{for }\alpha=1,2,\ldots,d,\end{equation}
and $x_t(\lambda)$ is on the cool part of the caustic.
\end{theorem}
\begin{proof}
At real turbulent times  there is a change in the cardinality,
\[\#\left\{
\lambda_d  =  \lambda_{d}(\lambda_1,\ldots,\lambda_{d-1}):
f_{(x_{t}(\lambda_1,\ldots,
\lambda_{d-1},\lambda_{d}(\lambda_1,\ldots,\lambda_{d-1})),t)}(\lambda_1)=c\right\}.\qedhere\]
\end{proof}

\subsection{White noise in $d$-orthogonal directions}
We now consider the Burgers fluid under the potential $V(x)=0$ and
 the noise
$\sum_{\alpha =1}^{d}\nabla k_{\alpha }({x})W_{\alpha }(t)$ where
$W_{\alpha }$ are $d$-independent Wiener processes and $k_{\alpha
}({x})=x_{\alpha }$ with ${x}=(x_{1},x_{2},\ldots,x_{d}).$ The
Burgers equation is then,
\begin{equation}\label{te1a}\frac{\partial v^{\mu}}{\partial t}
+(v^{\mu}\cdot\nabla)v^{\mu} = \hf[\mu^{2}]\Delta v^{\mu}
-\epsilon {\dot{W}}(t),\end{equation}
where  $W(t)=\left(W_1(t),W_2(t),\ldots,W_d(t)\right)$.
\begin{prop}\label{t3}
The stochastic action corresponding to the Burgers equation (\ref{te1a})
is,
\begin{eqnarray*}
\mathcal{A}({x_{0}},{x},t)
& = &
\frac{|x-x_0|^{2}}{2t}
+\frac{\epsilon}{t}(x-x_0)\cdot\int_{0}^{t}W(s)\di s
-\epsilon x\cdot W(t)
\\
&   & \quad
-\hf[\epsilon^{2}]\int_{0}^{t}|W(s)|^2\di s+\frac{\epsilon^{2}}{2t}\left|\int_{0}^{t}W(s)\di
u\right|^{2}
+S_{0}({x_{0}}).
\end{eqnarray*}
\end{prop}
\begin{proof}
The action is derived as in Section 2 using integration by parts.
\end{proof}

\begin{prop}\label{t6}
 If $x_t^{\epsilon}(\lambda)$ denotes the pre-parameter\-is\-ation of the random caustic for the stochastic Burgers equation (\ref{te1a}) and $x_t^0(\lambda)$ denotes the pre-parameter\-isation of the deterministic caustic (the $\epsilon=0$ case) then,
\[x_t^{\epsilon}(\lambda)=x_t^0(\lambda)-\epsilon\int_0^t W(u) \di u.\]
\end{prop}
\begin{proof}
Follows from Definition \ref{i2} and Theorem \ref{sflow}.
\end{proof}
Using Propositions  \ref{t3} and \ref{t6}, we can find the zeta process explicitly.
\begin{theorem}\label{t7}
In $d$-dimensions, the zeta process for the stochastic Burgers equation (\ref{te1a}) is,
\[
\zeta_t^c  = f_{(x_t^0(\lambda),t)}^0(\lambda_1) -\epsilon
x_t^0(\lambda)\cdot W(t) +\epsilon^2 W(t)\cdot \int_0^t W(s)\di s
-\hf[\epsilon^2]\int_0^t |W(s)|^2 \di s-c,\]
where $f_{(x,t)}^0(\lambda_1)$ is the deterministic reduced action
function, $x_t^0(\lambda)$ is the deterministic caustic and
$\lambda$ must satisfy the stochastic equation,
\begin{equation}\label{te5a}\nabla_{\lambda}\left(f^0_{(x_t^0(\lambda),t)}(\lambda_1)-
\epsilon x_t^0(\lambda)\cdot W(t)\right)
=0.\end{equation}
\end{theorem}
\begin{proof}
Follows from Theorem \ref{t2} having derived by induction the reduced action function from Theorem \ref{t3}.
\end{proof}
Equation (\ref{te5a}) shows that the value of $\lambda$ used in the zeta process may be either deterministic or random.
In the two dimensional case this gives,
\begin{equation}
0 =
\left( \nabla_x f^0_{(x_t^0(\lambda),t)}(\lambda_1)-
\epsilon  W(t)\right)\cdot\frac{\di x_t^0}{\di\lambda}(\lambda),\label{te5b}
\end{equation}
which  has a deterministic solution for $\lambda$ corresponding to a
cusp on the deterministic caustic. This point  will be returned to
in Section 5.3.

\subsection{Recurrence, Strassen and Spitzer}
One of the key properties associated with turbulence is the
intermittent recurrence of short intervals during which the fluid
velocity varies infinitely rapidly.

Using the law of the iterated logarithm, it is a simple matter to
show formally that if there is a time $\tau$ such that
$\zeta_{\tau}^c=0$, then there will be infinitely many zeros of
$\zeta_t^c$ in some neighbourhood of $\tau$. This will make the set
of zeros of $\zeta_t^c$ a perfect set and will result in a short
period during which the fluid velocity will vary infinitely rapidly.
However, this formal argument is not rigorous as it will not hold on
some set of times $t$ of measure zero \cite{Reynolds}.

The intermittent recurrence of turbulence will be demonstrated if we can show
that there is an unbounded increasing sequence of times at
which the zeta process is zero.

We begin by indicating the derivation of Strassen's form of the law
of the iterated logarithm from the theory of large deviations
\cite{Strook, Varadhan}. Consider a complete separable metric space
$X$ with a family of probability measures $\mathbb{P}_{\epsilon}$
defined on the Borel sigma algebra of $X$.

\begin{defn}
The family of probability measures $\mathbb{P}_{\epsilon}$ obeys the
large deviation principle with a rate function $I$ if there exists a lower semicontinuous
function $I:X\rightarrow [0,\infty]$ where:
\begin{enumerate}
\item for each $l\in\mathbb{R}$ the set $\{x: I(x)\le l\}$ is compact in $X$,
\item for each closed set $C\subset X$,
$\limsup\limits_{\epsilon\rightarrow 0} \epsilon \ln \mathbb{P}_{\epsilon}(C)\le -\inf\limits_{x\in C} I(x),$
\item for each open set $G\subset X$,
$\liminf\limits_{\epsilon\rightarrow 0} \epsilon \ln \mathbb{P}_{\epsilon}(G)\ge -\inf\limits_{x\in G} I(x).$
\end{enumerate}
\end{defn}
Let $X=C_0[0,1]$ where $C_0[0,1]$ is the space of continuous
functions $f:[0,1]\rightarrow \mathbb{R}^d$ with $f(0)=0$.  Let
$W(t)$ be a $d$-dimensional  Wiener process and
$\mathbb{P}_{\epsilon}$ be the distribution of $\sqrt{\epsilon}W(t)$
so that $\mathbb{P}_1$ is the Wiener measure.

\begin{theorem}
For the measure $\mathbb{P}_{\epsilon}$ the large deviation principle holds with a rate function,
\[I(f)=\left\{\begin{array}{lcl}
\hf\int_0^1 \dot{f}(t)^2\di t & : & f(t)\mbox{ absolutely continuous and } f(0)=0,\\
\infty & : & otherwise.
\end{array}\right.\]
\end{theorem}

\begin{defn}
The set of Strassen functions is defined by,
\[K = \left\{f\in C_0[0,1] :\quad 2I(f)\le 1\right\}.\]
\end{defn}
\begin{theorem}[Strassen's Law of the Iterated Logarithm]
Let
$Z_n(t) = \left(2n \ln\ln n\right)^{-\hf} W(nt)$
for $n\ge 2$ and $0\le t\le 1$ where $W(t)$ is a $d$-dimensional
Wiener process. For almost all paths $\omega$ the subset
$\left\{Z_n(t):  n=2,3,\ldots\right\}$
is relatively compact with limit set $K$.
\end{theorem}

Following the ideas of RTW, this theorem can be applied to the zeta process to demonstrate its recurrence.

\begin{cor}\label{t12}
There exists an unbounded increasing sequence of times $t_n$ for which $Y_{t_n}=0$, almost surely, where,
\[
Y_t =   W(t)\cdot \int_0^t W(s)\di s
-\hf\int_0^t |W(s)|^2 \di s,
\]
and $W(t)$ is a $d$-dimensional Wiener process.
\end{cor}
\begin{proof}If $h(n)=(2n \ln\ln n)^{-\hf}$
and $x(t)\in K$ then  there exists an increasing sequence $n_i$ such
that,
$Z_{n_i}(t) = h(n_i)W(n_it)\rightarrow x(t),$
as $i\rightarrow \infty.$

Consider  the behaviour of each term in $h(n_i)^2 n_i^{-1} Y_{t}$.
Firstly, by applying Lebesgue's dominated convergence theorem,
\begin{eqnarray}
h(n_i)^2 n_i^{-1} W(n_i)\cdot \int_0^{n_i} W(s)\di s
 &\rightarrow&
 x(1)\cdot\int_0^1x(r)\di r,\label{te6}
\end{eqnarray}
and,
\begin{eqnarray}
h(n_i)^2 n_i^{-1} \int_0^{n_i} |W(s)|^2\di s
 &\rightarrow&
 \int_0^1|x(r)|^2\di r,\label{te7}
\end{eqnarray}
as $i\rightarrow\infty$.

Now let $x(t)=(x_1(t),x_2(t),\ldots,x_d(t))$ where
$x_{\alpha}(t)=d^{-\hf}t $ for each $\alpha=1,2,\ldots,d$. Therefore, from equations (\ref{te6}) and (\ref{te7}), there is an increasing sequence of
times $t_i$ such that,
\[
h(t_i)^2t_i^{-1} Y_{t_i}  \rightarrow  \hf-\frac{1}{6} = \frac{1}{3},\]
as $i\rightarrow \infty$.

Alternatively, let
\[
x_{\alpha}(t)=\left\{
\begin{array}{ccc}
(d)^{-\hf}t & : & 0\le t \le \frac{1}{3},\\
 (d)^{-\hf}\left(\frac{2}{3}-t\right)& : & \frac{1}{3}\le t\le 1,
\end{array}\right.\]
for $\alpha = 1,2,\ldots d.$
Therefore, using equations (\ref{te6}) and (\ref{te7}) there is an
increasing sequence of times $\tau_i$ such that,
\[
h(\tau_i)^2\tau_i^{-1} Y_{\tau_i}  \rightarrow -\frac{1}{54}-\frac{1}{27} = -\frac{1}{18}.\]

Thus, the sequence  $t_i$ is an unbounded increasing infinite sequence of times at which $Y_t>0$, and the sequence $\tau_i$ is an unbounded increasing infinite sequence of times at which $Y_t<0$.
\end{proof}

\begin{cor}\label{t13}
If
$h(t)^2t^{-1}f^0_{(x_t^0(\lambda),t)}(\lambda_1)\rightarrow 0$ and $h(t)t^{-1}\sum\limits_{i=0}^d x_t^{0_i}(\lambda)\rightarrow 0,$
then the zeta process $\zeta_t^c$ is recurrent.
\end{cor}

A stronger condition on recurrence can be found in the two dimensional case if we work with small $\epsilon$ and neglect terms of order $\epsilon^2$ so that,
\begin{equation}\label{te7a}
\zeta_t^c  = f_{(x_t^0(\lambda),t)}^0(\lambda_1)
-\epsilon x_t^0(\lambda)\cdot W(t)
-c.
\end{equation}
 For this we use Spitzer's theorem, a proof of which can be found in Durrett \cite{Durrett}.
\begin{theorem}[Spitzer's Theorem]\label{t17} Let $D(t)=D_{1}(t)+iD_{2}(t)$ be a
complex Brownian motion where $D_{1}$ and $D_{2}$ are independent,
$D_{1}(0)=1$ and $D_{2}(0)=0$. Define the process $\theta_{t}$ as
the continuous process where $\theta_{0}=0$ and $\sin
(\theta_{t})=\frac{D_{2}(t)}{|D(t)|}$. Then, as
$t\rightarrow\infty$,
$$\mathbb{P}\left\{\frac{2\theta_t}{\ln t}\le
y\right\}\rightarrow\frac{1}{\pi}\int_{-\infty}^{y}\frac{\di x}{1+x^2}_.$$
\end{theorem}
The process $\theta_t$ gives the angle
swept out by $D(t)$ in time $t$, counting anti-clockwise loops as
$-2\pi$ and clockwise loops as $2\pi$.

Let
$A:\mathbb{R}^{+}\rightarrow\mathbb{R}^2$ and  consider the behaviour of the process, $Y_t = A(t)\cdot W(t).$
Assuming that $A(t)\neq 0$, let $\phi_t$ and $\theta_t$ measure the windings around the origin of
$A(t)$  and $W_t$ respectively. Then,
$Y_t = \epsilon|A(t)||W(t)|\cos
(\phi_{t}-\theta_{t}).$
 Therefore, for $Y_t=0$ we require  $\cos (\phi_{t}-\theta_{t})=0$, so that the
two vectors $A(t)$ and $W(t)$ are perpendicular to each
other. (Alternatively, this would be satisfied trivially if $A(t)$ were periodically zero with $t$.)

\begin{cor}\label{t22}
The small noise zeta process (\ref{te7a}) is recurrent if there exists a bounded function $h(t)$ where $h:\mathbb{R}^+\rightarrow\mathbb{R}^+$ such that,
$$h(t)\left(f_{(x_t^0(\lambda),t)}^0(\lambda_1)-c\right)\rightarrow 0,$$
as $t\rightarrow\infty$ and
there exists a function $n_t$ such that $n_t\rightarrow\infty$ with,
  $$\frac{4\pi^{2}n_{t}^{2}-\phi_{t}^{2}}{(\ln t )^{2}}
    <\frac{1}{4}\quad\mbox{and}\quad
    \frac{n_{t}\ln t}{16\pi^2n_t^2-4\phi_t^2+(\ln t)^{2}}
    \rightarrow 0,$$
as $t\rightarrow\infty$ where $A(t) = \epsilon h(t) x_t^0(\lambda)$.
\end{cor}

\subsection{Examples in two and three dimensions}
We now consider an explicit example in two dimensions. Since the parameter  $\lambda\in\mathbb{R}$, equation (\ref{te1b}) reduces to,
\begin{eqnarray}
0 &=& \frac{\di}{\di\lambda} f_{(x_t(\lambda),t)}(\lambda)
  =  \nabla_x f_{(x_t(\lambda),t)}(\lambda)\cdot
 \frac{\di x_t}{\di \lambda}(\lambda).\label{te13a}
 \end{eqnarray}
This gives three different forms of turbulence:
\begin{enumerate}
\item `zero speed turbulence' where
$\nabla f_{(x,t)}(x_0^1)=\dot{X}(t)=0.$ and so
 the Burgers fluid has zero velocity.

\item `orthogonal turbulence' where  $\nabla f_{(x_t(\lambda),t)}(\lambda)$  is orthogonal to
$\frac{\di x_t}{\di\lambda}(\lambda)$ so that the caustic tangent is orthogonal to the
Burgers fluid velocity.

\item `cusped turbulence' where
$\frac{\di x_t}{\di \lambda}(\lambda)=0,$
so there is a generalised cusp on the caustic at
$x_t(\lambda)$.
\end{enumerate}

As discussed previoiusly, cusped turbulence will occur at deterministic values of $\lambda$ and will also correspond to points of swallowtail perestroika on the level surfaces. As such, it is not only the simplest form to analyse, but also the most important.
The categorisation of turbulence leads to a
factorisation of equation (\ref{te13a}).
\begin{eg}[The generic Cusp]\label{t15}
For the generic Cusp, the zeta process reduces to,
\begin{eqnarray*}
\zeta_{t}^c& = &- \frac{3 \lambda^4 t}{8} + \frac{\lambda^6 t^3}{2}
- \epsilon\left(
 \lambda^3 t^2 W_{1}(t) - \frac{W_{2}(t)}{t} + \hf[3] \lambda^2 t
 W_{2}(t) \right)\\
&   &   +\epsilon^{2}\left(W(t)\cdot\int_0^tW(s) \di s
    -\hf\int_0^t|W(s)|^2\di s\right)-c,
\end{eqnarray*}
where $\lambda$ must be a root of,
$$0  = \hf[3]\lambda t (2\lambda^{4}t^{2}-\lambda^{2}
-2\epsilon\{\lambda tW_{1}(t)+W_{2}(t)\}).$$ The factor $\lambda=0$
corresponds to the cusp on the caustic while the roots of the second
factor correspond to orthogonal and zero speed turbulence. Firstly,
if $\lambda=0$ then from Corollary \ref{t13}, the zeta process is
recurrent. Therefore, the turbulence occurring at the cusp on the
generic Cusp caustic is recurrent.

Alternatively, for large times it can be shown formally that the
four roots which give rise to orthogonal and zero speed turbulence
all tend towards zero. Thus all four roots tend towards the cusp and
consequently, the zeta processes associated with each root will be
recurrent.

Moreover, it can be shown that the whole caustic is cool and so all of these points of turbulence will be genuine.\end{eg}

Next consider the three dimensional case. Thus
$\lambda\in\mathbb{R}^2$ and equation (\ref{te1b}) becomes the pair,
\begin{equation}
 0 =  \nabla_x f_{(x_t(\lambda),t)}(\lambda_1)\cdot
 \frac{\di x_t}{\di \lambda_1}(\lambda),\quad
 0 =
 \nabla_x f_{(x_t(\lambda),t)}(\lambda_1)\cdot
 \frac{\di x_t}{\di \lambda_2}(\lambda).\label{te15b}
 \end{equation}

In direct correlation to the two dimensional case, we can categorise
three dimensional turbulence depending on how we solve equations
(\ref{te15b}):
\begin{enumerate}
\item `zero speed turbulence' where again $\nabla_x f_{(x_t(\lambda),t)}(\lambda_1) = 0$,
\item `orthogonal turbulence' where all three vectors $\nabla_x f_{(x_t(\lambda),t)}(\lambda_1)$, $\frac{\di x_t}{\di \lambda_1}(\lambda)$ and $\frac{\di x_t}{\di \lambda_2}(\lambda)$ are mutually orthogonal.
\item `subcaustic turbulence' where the vectors $\frac{\di x_t}{\di \lambda_1}(\lambda)$ and $\frac{\di x_t}{\di \lambda_2}(\lambda)$ are linearly dependent. The ``\emph{subcaustic}" is the region of the caustic where the tangent space drops one or more dimensions. In three dimensions it corresponds to folds in the caustic.
\end{enumerate}
As in the two dimensional case, it follows from Proposition \ref{t6}
that the values of $\lambda$ that determine the subcaustic are
deterministic.  However, unlike the two dimensional case, subcaustic
turbulence only occurs at points where the Burgers fluid velocity is
orthogonal to the subcaustic. Hence, we are selecting random points
on a deterministic curve, and so subcaustic turbulence involves
random values of $\lambda$. Again the categoristaion of turbulence
leads to a factorisation in equations (\ref{te15b}).
\begin{eg}[The butterfly]
Let $S_0(x_0,y_0,z_0)=x_0^3y_0+x_0^2z_0$, this gives a butterfly caustic -- the three dimensional analogue of the generic Cusp.
The zeta process is,
\begin{eqnarray*}
\zeta_t^c & = &
    \lambda_1^3\lambda_2 - \frac{3}{2}\lambda_1^4t - 4\lambda_1^6t + \frac{9}{2}\lambda_1^4\lambda_2^2t - 12\lambda_1^5\lambda_2t^2 -
  27\lambda_1^7\lambda_2t^2 + 8\lambda_1^6t^3  \\
 &  & + 36\lambda_1^8t^3+ \frac{81}{2}\lambda_1^{10}t^3+
  \epsilon \Big(\big( 3\lambda_1^2\lambda_2t - 4\lambda_1^3t^2 - 9\lambda_1^5t^2 \big)
      W_1(t)  \\
      &   &\qquad- \big(\lambda_2 + \lambda_1^3t \big) W_2(t)+
     \big( 3\lambda_1\lambda_2 + \frac{1}{2t} - 3\lambda_1^2t - \frac{9}{2}\lambda_1^4t \big)
      W_3(t) \Big)\\
      &   &+\epsilon^{2}
\Big(W(t)\cdot\int_0^tW(s) \di s
        -\hf\int_0^t|W(s)|^2\di s\Big)-c,
\end{eqnarray*}
where  $\lambda=(\lambda_1,\lambda_2)$ must satisfy,
\begin{eqnarray*}
    0 & = &
       135t^{3}\lambda_1^{9}+96t^{3}\lambda_1^{7}-63\lambda_2t^{2}\lambda_1^{6}
    +\big(16t^{3}-8t\big)\lambda_1^{5}-\left(20\lambda_2+15\epsilon
    W_{1}(t)\right)t^{2}\lambda_1^{4}\\
    &  &\quad+
    \left(6\lambda_2^{2}-2-6\epsilon W_{3}(t)\right)t \lambda_1^{3}+
    \left(\lambda_2-4t^{2}\epsilon W_{1}(t)-t\epsilon
    W_{2}(t)\right)\lambda_1^{2}\\
    &   &\quad
    +2\left(\lambda_2 W_{1}(t)-W_{3}(t)\right)\epsilon t \lambda_1
        +\epsilon \lambda_2 W_{3}(t),\\
   0 & = &
    -27t^{2}\lambda_1^{7}-12t^{2}\lambda_1^{5}+9t \lambda_1^{4}\lambda_2+\lambda_1^{3}+3t\epsilon
    W_{1}(t)\lambda_1^{2}+3\epsilon W_{3}(t)\lambda_1-\epsilon
    W_{2}(t).
\end{eqnarray*}
Eliminating $\lambda_2$ gives the factorisation,
\begin{eqnarray}
0& =
&\left(54t^{2}\lambda_1^{7}+6t^{2}\lambda_1^{5}+\lambda_1^{3}+3t\epsilon
W_{1}(t) \lambda_1^{2} +3\epsilon W_{3}(t)\lambda_1-\epsilon
W_{2}(t) \right)
\nonumber\\
&   &\times\left(\lambda_1^{3}-3\epsilon W_{3}(t)\lambda_1+2\epsilon
W_{2}(t)\right), \label{te22}
\end{eqnarray}
where the first factor gives zero speed and orthogonal turbulence while the second factor gives subcaustic turbulence.

For large times, it can be shown formally that of the seven roots
corresponding to zero and orthogonal turbulence, five should tend to
$\lambda=(0,0)$ and so should give a recurrent zeta process. None of
the remaining roots give recurrence \cite{Neate2}.
\end{eg}
\subsection{The harmonic oscillator potential}
It is not always necessary to resort to Strassen's law to
demonstrate the recurrence of turbulence; some systems have an
inherent periodicity which produces such behaviour. The following
two dimensional example is taken from RTW \cite{Truman6} in which a
single Wiener process acts in the $x$ direction.
\begin{eg}
Let $k_t(x,y) = x$, $V(x,y) = \hf (x^2\omega_1^2+ y^2\omega_2^2)$
and $S_0(x_0,y_0) = f(x_0)+g(x_0)y_0$ where $f$, $f'$, $f'''$, $g$,
$g'$, $g'''$ are zero when $x_0=\alpha$ and $g''(\alpha) \neq 0$.
Then the zeta process for turbulence at $\alpha$ is given by,
\begin{eqnarray*}
\zeta_t^c & = &
-\frac{\omega_2}{4 g''(\alpha)}\sin(2\omega_2 t)
\csc^2(\omega_1 t)\left\{\sin (\omega_1 t)f''(\alpha)+\omega_1\cos(\omega_1 t)\right\}^2\\
&   & \quad +\epsilon \csc(\omega_1 t)R_t -\frac{1}{4}\alpha^2\omega_1\sin(2\omega_1 t)-c,
\end{eqnarray*}
where $R_t$ is a stochastic process which is well defined for all $t$.

Therefore, $\zeta_t\rightarrow\pm\infty$ as $t\rightarrow
\frac{k\pi}{\omega_1}$ because $\csc^2 (k\pi)=\infty$ where the sign
depends upon the sign of $-\frac{\sin(2\omega_2 t)}{g''(\alpha)}$.
Thus, it is possible to construct an unbounded increasing sequence
of times at which $\zeta_t$ switches between $\pm\infty$ and so by
continuity and the intermediate value theorem there will almost
surely exist an increasing unbounded sequence $\{t_k\}$ at which
$\zeta_{t_k}=0$.
\end{eg}
\section{Complex turbulence}

We now consider a completely different approach to turbulence based
on the work of Section 4.  Let
$\left(\lambda,x^2_{0,\mathrm{C}}(\lambda)\right)$ denote the
parameterisation of the pre-caustic so that
$x_t(\lambda)=\Phi_t\left(\lambda,x^2_{0,\mathrm{C}}(\lambda)\right)$
is the pre-parameterisation of the caustic. When,
\[Z_t=\mbox{Im}\left\{\Phi_t(a+\rmi\eta, x^2_{0,\mathrm{C}}(a+\rmi\eta))\right\},\] is random,
 the values of $\eta(t)$ for which $Z_t=0$ will form a stochastic process.
 The zeros of this new process will correspond to points at which the real pre-caustic touches the complex pre-caustic.
 \begin{defn}
 The complex turbulent times $t$ are defined to be times $t$ when the real and complex pre-caustics touch.
 \end{defn}

  The points at which these surfaces touch correspond to swallowtail perestroikas on the caustic.
   \begin{theorem} \label{s7}  Let
    $x_{t}(\lambda)$ denote the pre-parameterisation of the caustic \linebreak[4]
    where
    $\lambda\in\mathbb{R}$ and $x_t(\lambda)$ is a real analytic function. If at time
$\tilde{t}$ a swallowtail perestroika occurs on the caustic when $\lambda=\tilde{\lambda}$ then,
\[f'_{(x_{\tilde{t}}(\tilde{\lambda}),{\tilde{t}})}(\tilde{\lambda})=
f''_{(x_{\tilde{t}}(\tilde{\lambda}),{\tilde{t}})}(\tilde{\lambda})=
f'''_{(x_{\tilde{t}}(\tilde{\lambda}),{\tilde{t}})}(\tilde{\lambda})=
f^{(4)}_{(x_{\tilde{t}}(\tilde{\lambda}),{\tilde{t}})}(\tilde{\lambda})=0.\]
\end{theorem}

Assuming that $f_{(x,t)}(x_0^1)$ is a polynomial in $x_0^1$ we can
use the resultant to state explicit conditions for which this holds.
\begin{lemma}
Let $g$ and $h$ be polynomials of degrees $m$ and $n$ respectively with no common roots or zeros.
Let $f=gh$ be the product polynomial. Then the resultant,
\[ R(f,f')=(-1)^{mn} \left(\frac{m!n!}{N!}\frac{f^{(N)}(0)}{g^{(m)}(0)h^{(n)}(0)}\right)^{N-1}
R(g,g') R(h,h') R(g,h)^2,\]
where $N=m+n$ and $R(g,h)\neq 0$.
\end{lemma}
\begin{proof}
See \cite{Neate}.
\end{proof}

 Since $f'_{(x_t(\lambda),t)}(x_0^1)$ is a polynomial in $x_0$
with real coefficients, its zeros are real or occur in complex
conjugate pairs. Of the real roots, $x_0=\lambda$ is repeated. So,
\[f'_{(x_t(\lambda),t)}(x_0^1) = (x_0^1-\lambda)^2
Q_{(\lambda,t)}(x_0^1)
 H_{(\lambda,t)}(x_0^1),\]
 where $Q$ is the product of quadratic factors,
\[Q_{(\lambda,t)}(x_0^1)=\prod\limits_{i=1}^{q}\left\{(x_0^1-a_t^i)^2+(\eta_t^i)^2\right\},\]
 and
 $H_{(\lambda,t)}(x_0^1)$ the product of real factors corresponding to real zeros.
This gives,
\[\left.f'''_{(x_t(\lambda),t)}(x^1_0)\right|_{x^1_0=\lambda}
= 2\prod\limits_{i=1}^q \left\{(\lambda-a_t^i)^2+(\eta_t^i)^2\right\}H_{(\lambda,t)}(\lambda).\]

We now assume that the real roots of $H$ are distinct as are the complex roots of $Q$. Denoting
$\left.f'''_{(x_t(\lambda),t)}(x_0^1)\right|_{x_0^1=\lambda}$ by $f'''_t(\lambda)$ etc, a simple calculation gives
\begin{eqnarray*}
\lefteqn{\left|R_{\lambda}(f'''_t(\lambda),f^{(4)}_t(\lambda))\right|=}\\
 &  &\!\!\!\!\!\!\! K_t \prod\limits_{k=1}^q(\eta_t^k)^2\prod\limits_{j\neq k}
\left\{(a_t^k-a_t^j)^4+2((\eta_t^k)^2+(\eta_t^j)^2)(a_t^k-a_t^j)^2+ ((\eta_t^k)^2-(\eta_t^j)^2)^2\right\}\\
& &\quad\times \left|R_{\lambda}(H,H')\right|\left|R_\lambda(Q,H)\right|^2,
\end{eqnarray*}
$K_t$ being a positive constant. Thus, the condition for a swallowtail perestroika to occur is that
\[\rho_{\eta}(t):=\left|R_{\lambda}(f'''_t(\lambda),f^{(4)}_t(\lambda))\right|=0,\]
where we call $\rho_{\eta}(t)$ the \emph{resultant eta process}.

When the  zeros of $\rho_{\eta}(t)$ form a perfect set, swallowtails
will spontaneously appear and disappear on the caustic infinitely
rapidly. As they do so, the geometry of the caustic will rapidly
change Moreover, Maxwell sets will be created and destroyed with
each swallowtail that forms and vanishes as when a swallowtail forms
it contains a region with two more pre-images than the surrounding
space.  This will add to the turbulent
 nature of the solution in these regions.
We call this `complex turbulence' occurring at the turbulent times
 which are the zeros of the resultant eta process.

Complex turbulence can be seen as a special case of real turbulence
which occurs at specific generalised cusps of the caustic. Recall
that when a swallowtail perestroika occurs on a curve, it also
satisfies the conditions for having a generalised cusp. Thus, the
zeros of the resultant eta process must coincide with some of the
zeros of the zeta process for certain forms of cusped turbulence. At
points where the complex and real pre-caustic touch, the real
pre-caustic and pre-level surface touch in a particular manner (a
double touch) since at such a point two swallowtail perestroikas on
the level surface have coalesced.

Thus, our separation  of complex turbulence from real turbulence can
be seen as an alternative form of categorisation to that outlined in
Section 5.3  which could be extended to include other perestroikas.

\section*{Acknowledgement}
One of us (AT) would like to record his indebtedness to John T Lewis
as his teacher, mentor and friend. This paper could not have been
written without John's inspirational work on large deviations which
underlies our work.


\begin{thebibliography}{99}

\bibitem{Arnold4} Arnol'd V I, Shandarin S F and  Zeldovich Y B 1982
The large scale structure of the universe 1
{\it Geophys. Astrophys. Fluid Dyn.} {\bf 20} 111--30

\bibitem{Arnold3}  Arnol'd V I 1986 {\it Catastrophe Theory}
(Berlin: Springer-Verlag)

\bibitem{Arnold}  Arnol'd V I 1989 {\it Mathematical Methods of Classical Mechanics}
(New York: Springer-Verlag)

\bibitem{Arnold2}  Arnol'd V I 1990 {\it Singularities of Caustics and Wave
Fronts. Mathematics and its Applications (Soviet Series) 62}
(Dordrecht: Kluwer Academic Publishers Group)

\bibitem{Dafermos} Dafermos C 2000 {\it Hyperbolic Conservation
Laws in Continuum Physics. Grundlehren der Mathematischen
Wissenschaten 325} (Berlin: Springer-Verlag)

\bibitem{Davies2} Davies I M,  Truman A and Zhao H  2002
Stochastic heat and Burgers equations and their singularities I -
geometric properties
{\it J. Math. Phys.} \textbf{43} 3293-328

\bibitem{Davies3} Davies I M,  Truman A and Zhao H  2005 Stochastic heat
and Burgers equations and their singularities II. Analytical properties and limiting distributions {\it J.  Math. Phys.} \textbf{46} 043515

\bibitem{Durrett}
Durrett R 1984 {\it Brownian motion and martingales in analysis} (Belmont: Wadsworth)


\bibitem{E} E Weinan, Khanin K,  Mazel A and  Sinai Y 2000
Invariant measures for Burgers equations with stochastic forcing
{\it Ann. Math.} \textbf{151} 877-960

\bibitem{Elworthy} Elworthy K D, Truman  A and Zhao H
Stochastic elementary formulae on caustics 1: One dimensional linear
heat equations {\it UWS MRRS Preprint}

\bibitem{Elworthy2} Elworthy K D, Truman  A and Zhao H 2005
Generalised Ito formulae and space-time Lebesgue-Stieltjes integrals
of local times {\it To appear in Seminaire de Probabilites
Strasbourg Vol. 40}

  \bibitem{Freidlin} Freidlin M I and Wentzell A D 1998
{\it Random Perturbations of Dynamical Systems} (New York:
Springer-Verlag)

\bibitem{Gilmore} Gilmore R 1981 {\it Catastrophe Theory for
Scientists and Engineers} (New York: John Wiley)

\bibitem{Hwa} Hwa R C and Teplitz V L 1966
{\it Homology and Feynman integrals} (New York: W A Benjamin)

\bibitem{Kac} Kac M 1959 {\it Probability and Related Topics in Physical Science}
 (New York: Interscience Publishers)

\bibitem{Klein} Klein F 1922 {\it \"{U}ber den Verlauf der Abelschen
Integrale bei den Kurven vierten Grades} in {\it Gesammelte
Mathematische Abhandlungen II} ed Fricke R and Vermeil H (Berlin:
Springer)


\bibitem{Kolokoltsov} Kolokoltsov V N, Schilling R L and
Tyukov A E 2004 Estimates for multiple stochastic integrals and
stochastic Hamilton-Jacobi equations {\it Rev. Mat. Iberoamericana}
\textbf{20} 333-80



\bibitem{Kunita}  Kunita H 1984 {\it Stochastic differential
equations and stochastic flows of homeomorphisms}, in {\it
Stochastic Analysis and Applications, Advances in Probability and
Related Topics. Vol. 7} ed Pinsky M A (New York: Marcel Dekker)


\bibitem{Maslov} Maslov V P 1972 {\it Perturbation Theory and Asymptotic
Methods} (Paris: Dunod)

\bibitem{Maslov2} Maslov V P and  Fedoriuk M V 1981 {\it Semi-Classical
Approximation in Quantum Mechanics. Mathematical Physics and Applied
Mathematics Vol. 7} (Dordrecht: Riedel Publishing Company)

\bibitem{Neate2} Neate A D 2005 {\it A one dimensional analysis of the singularities of the $d$-dimensional stochastic Burgers equation} PhD thesis UWS

\bibitem{Neate} Neate A D and Truman A 2005 A one dimensional analysis of real and complex turbulence and the Maxwell set for the stochastic Burgers equation {\it J. Phys. A: Math. Gen.} \textbf{38} 7093--127

\bibitem{Neate3} Neate A D and Truman A 2005 The Maxwell set for the stochastic Burgers equation {\it In preparation.}

\bibitem{Reynolds} Reynolds C 2002 {\it On the polynomial
swallowtail and cusp singularities of stochastic Burgers equations}
PhD thesis UWS

\bibitem{Salmon} Salmon G 1934 {\it A Treatise on the Higher Plane Curves}
(New York: G E Stechert  Co)


\bibitem{Shandarin} Shandarin S F and Zeldovich Y B 1989
The large scale structure of the universe 2: turbulence,
intermittency, structures in a self gravitating medium
{\it Rev. Mod. Phys.} \textbf{6} 185-220

\bibitem{Strook}
Stroock D W 1984 {\it An introduction to the theory of large deviations} (New York: Springer-Verlag)

\bibitem{Truman} Truman A and Zhao H  1995 The stochastic Hamilton-Jacobi equations and related topics: a survey {\it LMS Lecture Note Ser. 216}
(Cambridge: Cambridge University Press) p~287

\bibitem{Truman2} Truman A and
Zhao H  1996 The stocahstic Hamilton-Jacobi equations, stochastic
heat equations and Schr\"{o}dinger equations {\it Stochastic
Analysis and Applications. Proc. of the 5th Gregynog Symp. held in
Powys July 9--14 1995 } ed. Davies I M et al (River Edge NJ: World
Scientific) p~441--64

 \bibitem{Truman3} Truman A and Zhao H  1996 On stochastic diffusion equations and stochastic Burgers equations {\it J.
Math. Phys.} \textbf{37} 283--307

\bibitem{Truman4} Truman A and Zhao H  1996 Quantum mechanics of
charged particles in random electromagnetic fields {\it J. Math.
Phys.} \textbf{37} 3180--97


\bibitem{Truman5} Truman A and Zhao H  1998
Stochastic Burgers equations and their semi classical expansions
{\it Comm. Math. Phys.} \textbf{194} 231-48



\bibitem{Truman6} Truman A,  Reynolds C N and Williams D 2003
Stochastic Burgers equations in $d$-dimensions -- a one dimensional
analysis: Hot and cool caustics and intermittence of stochastic
turbulence {\it Probabilistic Methods in Fluids} ed Davies I M et al
(Singapore: World Scientific) pp~239--62

\bibitem{vander}Van Der Waerden 1949 {\it Modern Algebra} Vols. 1 and 2 (New York: Frederick Ungar Publishing)
\bibitem{Varadhan} Varadhan S R S 1984 {\it Large deviations and applications} (Philadelphia: SIAM)

\end{thebibliography}
\end{document}